\documentclass[%
 aip,
 amsmath,amssymb
 reprint,%
]{revtex4-1}
\usepackage{graphicx}
\usepackage{dcolumn}
\usepackage{bm}

\usepackage[dvipsnames]{xcolor}
\usepackage{soul}
\usepackage[utf8]{inputenc}
\usepackage[T1]{fontenc}
\usepackage{mathptmx}
\usepackage{etoolbox}


\begin{document}


\title{Heating from Above in Non-scattering Suspensions: Phototactic Bioconvection under Collimated Irradiation}
\author{Sandeep Kumar}
\thanks{Corresponding author}
\email{sandeepkumar1.00123@gmail.com}
\affiliation{
Department of Mathematics, PDPM  Indian Institute of Information Technology Design and Manufacturing, Jabalpur 482005, India
}%

\author{Shaowei Wang}
\email{shaoweiwang@sdu.edu.cn}
\affiliation{
Department of Engineering Mechanics, School of Civil Engineering, Shandong University, Jinan 250061, PR China
}


\date{\today}

\begin{abstract}

Examining phototactic bioconvection in non-scattering suspensions with upper heating and collimated irradiation, this study delves into the intricate dynamics influenced by light and microorganisms. The study focuses on the linear stability of the basic state, examining neutral curves. The numerical analysis involves solving a system of equations using the MATLAB bvp4c solver. The investigation considers the impact of parameters, such as the thermal Rayleigh number, critical total intensity, and Lewis number, on the critical bioconvection Rayleigh number. As the critical total intensity varies, a transition from a stationary to an oscillatory solution (and vice versa) is observed. Phototactic microorganisms are incorporated into the model, and results show how varying parameters affect convection patterns and stability. The findings reveal interesting phenomena, including Hopf bifurcations and limit cycles.

\end{abstract}

\maketitle

\section{Introduction}
The phenomenon of bioconvection, characterized by the spontaneous formation of patterns in suspensions of swimming microorganisms, has been a subject of extensive research since its initial exploration by Wager \cite{ref2}. Platt \cite{ref3} coined the term ``bioconvection'', highlighting the non-predetermined nature of the patterns generated by these microorganisms. The aggregation of slightly denser microorganisms due to tactic behavior, coupled with their upward movement, leads to the emergence of distinct patterns, disappearing when the microorganisms cease swimming \cite{ref1}.

Various microorganisms, such as flagellated green algae (\textit{Volvox}, \textit{Euglena}, \textit{Dunaliella}, and \textit{Chlamydomonas}), exhibit bioconvection patterns influenced by factors like depth, microorganism concentration, and motion \cite{ref6,ref5,ref4,ref2}. The responses of microorganisms to external stimuli, collectively termed taxes, encompass phenomena like chemotaxis, gravitaxis, gyrotaxis, and phototaxis. Phototaxis, explored in this work, involves motion either toward or away from light. Experiments shed light on how light intensity and gradient influence bioconvection patterns \cite{ref5,ref2}, with pattern characteristics changing based on the amount of light \cite{ref7}. The interplay between total light intensity $\mathcal{G}$, a critical value $\mathcal{G}_c$, and light absorption by microorganisms contributes to pattern variations \cite{ref8}.

In suspensions illuminated from above with finite depth, a fundamental equilibrium state arises, marked by a balance between phototaxis induced by diffusion and light absorption. This equilibrium manifests as a horizontal sublayer densely populated by microorganisms, with the region above remaining gravitationally stable and the region below undergoing gravitational instability. The critical total intensity $\mathcal{G}_c$ determines the sublayer's location: above, below, or in the middle of the suspension, impacting the occurrence of penetrative convection in unstable fluid layers \cite{ref10}.

While early studies primarily focused on constant-temperature fluid suspensions, the present understanding extends to biothermal convection. Researchers like Kuznetsov \cite{kuznetsov2005onset, kuznetsov2011non}, Zhao et al. \cite{zhao2018linear, zhao2019darcy}, Balla et al. \cite{balla2020bioconvection}, and Hussain et al. \cite{hussain2022thermal} delved into the instability of gyrotactic microorganisms in heated suspensions, exploring both stationary and oscillatory convection. The inclusion of nanoparticles, porous media, and thermal radiation expanded the scope of these studies. The introduction of phototactic bioconvection was pioneered by Vincent and Hill \cite{ref9}, focusing on non-scattering absorbing suspensions. Subsequent models incorporated isotropic scattering \cite{ref14}, isotropic scattering with rigid boundaries \cite{kumar2023}, forward anisotropic scattering \cite{ref15}, and both diffuse and collimated light \cite{ref16,panda2020effects}. Investigations explored the impact of oblique light \cite{ref17} and scattering albedo variations \cite{ref18}, with recent work considering the effect of rotation on non-scattering medium \cite{kumar2023effect}. However, no theoretical research has comprehensively addressed both phototactic and temperature gradient processes with heating from above.

In this study, we adopt the phototaxis model by Vincent and Hill \cite{ref9}, utilizing the Navier-Stokes equation and a cell conservation equation to model an incompressible fluid within a phototactic suspension heating from above or cooling from below with illumination from above. With motile algae exhibiting phototaxis behavior and being reliant on photosynthesis for nutrition, our investigation aims to realistically capture the interplay between phototaxis and thermal influences.

The structure of the paper unfolds as follows: Section \ref{sec2} presents the mathematical model for phototactic bioconvection and corresponding boundary conditions. Section \ref{scaling} presents scaled governing equations. The steady-state solution is discussed in Section \ref{sec3}, followed by a linear stability analysis and normal mode analysis in Section \ref{sec4}-\ref{normal}. Numerical results and conclusions are presented in Sections \ref{sec5} and \ref{sec6}.

\begin{figure*}
    \centering
    \includegraphics[width=16cm, height=10cm ]{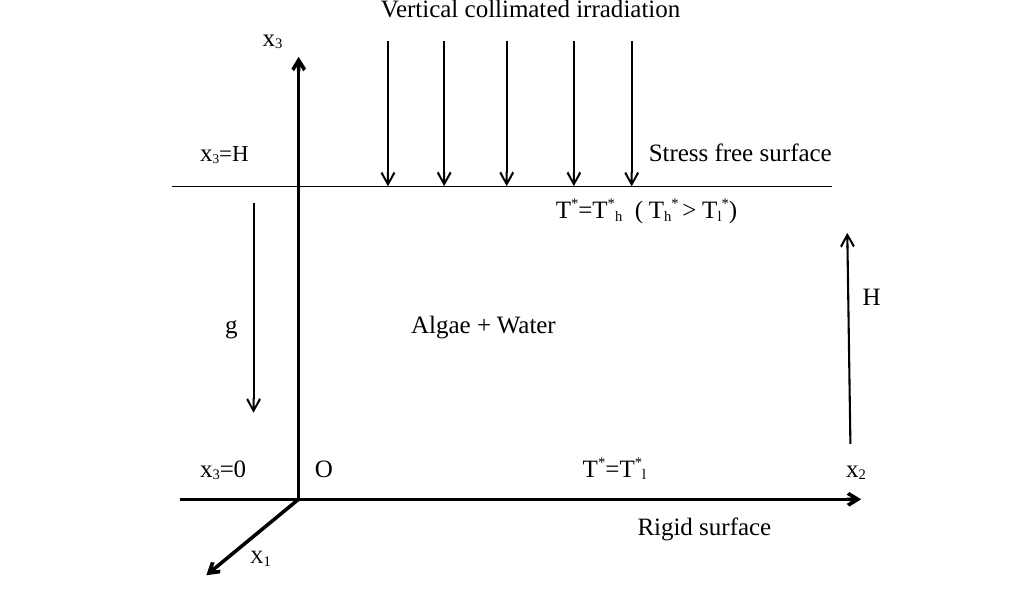}
    \caption{The geometry of the problem. }
   \label{fig:thermalfig.eps}
 \end{figure*}

\section{MATHEMATICAL MODEL}
\label{sec2}
Consider a Cartesian coordinate system in three dimensions denoted as \(O x_1 x_2 x_3\), extending infinitely in width but with a finite depth \(H\). The lower boundary at \(x_3=0\) is treated as rigid, while the upper boundary at \(x_3=H\) is stress-free. Direct collimated lighting illuminates the top boundary, with no consideration for diffused or oblique irradiation. Rotation's impact on the suspension is also neglected. The intensity at a point \(\boldsymbol{x}=(x_1,x_2,x_3)\) in a particular direction \(\boldsymbol{r}\) is described by the function \(I(\boldsymbol{x},\boldsymbol{r})\), with components \(\sin{\alpha_1}\cos{\alpha_2}\) in the \(\boldsymbol{i}\) direction, \(\sin{\alpha_1}\sin{\alpha_2}\) in the \(\boldsymbol{j}\) direction, and \(\cos{\alpha_1}\) in the \(\boldsymbol{k}\) direction. Here, \(\alpha_1\) is the polar angle, and \(\alpha_2\) is the azimuthal angle.

The continuity equation is given by
\begin{equation}
\label{eqn:equation8}
div(\boldsymbol{v})=0.
\end{equation}

The momentum equation, related to Boussinesq, is expressed as \cite{zhao2018linear}
\begin{align}
\label{eqn:equation9}
\varrho\left(\frac{\partial}{\partial t}+\boldsymbol{v}\cdot\nabla\right)\boldsymbol{v} &= \mu\nabla^2\boldsymbol{v} - \boldsymbol{\nabla} \mathcal{P} - n\vartheta \Delta \varrho g \boldsymbol{k}  - \varrho (1-\beta(T^*-T_h^*))g \boldsymbol{k}.
\end{align}

Here, \(t\) is time, \(\boldsymbol{v}\) is fluid velocity, \(\mu\) is dynamic viscosity, \(\rho\) is fluid density, \(\Delta \rho\) is the density of each algal cell (\(\Delta \rho/\rho \ll 1\)), \(\vartheta\) is the volume of each algal cell, \(g\) is gravitational acceleration, \(n\) is cell concentration, \(\mathcal{P}\) is excess pressure above hydrostatic, \(\beta\) is the volumetric thermal expansion coefficient, \(T^*\) is the temperature of the fluid, \(T_h^*\) is the temperature of the upper wall, and \(T_l^*\) is the temperature of the lower wall.

The thermal energy equation is given by
\begin{equation}
\varrho c \left[\frac{\partial T^*}{\partial t}+ \boldsymbol{v}\cdot\boldsymbol{\nabla} T^*\right] = \boldsymbol{\nabla}(k \boldsymbol{\nabla} T^*).
\end{equation}

The conservation of cells equation is given by \cite{ref1,ref9}:
\begin{eqnarray}
\label{eqn:equation10}
\frac{\partial n}{\partial t} = -\boldsymbol{\nabla}\cdot\boldsymbol{A}, \nonumber\\ 
\boldsymbol{A} = n \boldsymbol{U}_c + n\boldsymbol{v} - \boldsymbol{D}\cdot\boldsymbol{\nabla} n,
\end{eqnarray}
where \(\boldsymbol{A}\) represents the flux of cells, including the mean cell swimming velocity \(\boldsymbol{U}_c\), the flux caused by advection, and the random component of cell locomotion. \(\boldsymbol{D}=D\boldsymbol{I}\) is the isotropic and constant diffusivity tensor, where \(D\) is the diffusion coefficient and $\boldsymbol{I}$ is the identity tensor.

The mean cell swimming velocity is given by
\[
\label{eqn:equation5}
\boldsymbol{U}_c = U_c\bar{\boldsymbol{P}},
\]
where \(U_c\) is the average cell swimming speed, and \(\bar{\boldsymbol{P}}\) is the average cell swimming direction, which is calculated as \cite{ref9}
\begin{equation}
\label{eqn:equation6}
    \bar{\boldsymbol{P}}=T(\mathcal{G})\boldsymbol{k}.
\end{equation}

The taxis function \(T(\mathcal{G})\) is defined as:
\[
\label{eqn:equation7}
T(\mathcal{G}) \left\{ \begin{array}{lll}
\ge 0 & \text{for $\mathcal{G}_c \ge \mathcal{G}$};\\
< 0 & \text{for $\mathcal{G}_c < \mathcal{G}$}.
\end{array} \right. 
\]

The precise expression of the phototaxis function $T(\mathcal{G})$, varies among microorganisms due to their distinctive characteristics \cite{ref9}. To provide a further level of understanding of this concept, a particular implementation of the phototaxis function might be described as follows: \cite{ref13,ref17}
\[
T(\mathcal{G})=0.8\sin{\left(\frac{3}{2}\pi\varphi(\mathcal{G})\right)}-0.1\sin{\left(\frac{1}{2}\pi\varphi(\mathcal{G})\right)},
\]
where \(\varphi=\mathcal{G}e^{\chi(\mathcal{G}-1)}\), and \(-1.1\le \chi\le 1.1\) implies \(0.3\le\mathcal{G}_c\le 0.8\) \cite{ref17}.

To maintain mass and cell conservation in an incompressible fluid suspension, horizontal boundaries are required. The boundary conditions for the stress-free top surface are given by Eq. (\ref{eqn:equation12}), and for the rigid bottom surface:
\begin{align}
\label{eqn:equation12}
\boldsymbol{v}\cdot\boldsymbol{k} &= 0, \quad \frac{\partial^2}{\partial x_3^2}(\boldsymbol{v}\cdot\boldsymbol{k})=0, \quad T^*=T_h^*, \quad \boldsymbol{A}\cdot\boldsymbol{k}=0,  \quad \text{at} \quad x_3=H,
\end{align}
while for the rigid bottom surface 
\begin{equation}
\label{eqn:equation11}
\boldsymbol{v}=0, \quad \boldsymbol{v}\times\boldsymbol{k}=0, \quad  T^*=T_l^*, \quad\boldsymbol{A}\cdot\boldsymbol{k}=0,\quad \text{at} \quad x_3=0.
\end{equation}

The radiation transfer equation is given by \cite{ref-modest,ref-chand}
\begin{equation}
\label{eqn:equation3}
\boldsymbol{r}\cdot \boldsymbol{\nabla}I(\boldsymbol{x},\boldsymbol{r})+\psi I(\boldsymbol{x},\boldsymbol{r})=0,
\end{equation}
where the absorption coefficient \(\psi=\iota n\) is defined with the concentration \(n\). The boundary conditions exclude light reflection on the top and bottom surfaces, given as:
\begin{equation}
\label{eqn:equation4}
I(x_1,x_2,H,\alpha_1,\alpha_2)=I^0\delta(\boldsymbol{r}-\boldsymbol{r}^0) , \quad \pi/2\le\alpha_1\le\pi,
\end{equation}
and
\begin{equation}
I(x_1,x_2,0,\alpha_1,\alpha_2)=0,\quad 0\le\alpha_1\le\pi/2.
\end{equation}

The radiative heat flux \(\boldsymbol{q}(\boldsymbol{x})\) and total intensity \(\mathcal{G}(\boldsymbol{x})\) at a point \(\boldsymbol{x}=(x_1,x_2,x_3)\) are given by \cite{ref-modest}:
\begin{equation}
\label{eqn:equation2}
\boldsymbol{q}(\boldsymbol{x})=\int_0^{4\pi}I(\boldsymbol{x},\boldsymbol{r})\boldsymbol{r}d\omega,
\end{equation}
\begin{equation}
\label{eqn:equation1}
\mathcal{G}(\boldsymbol{x})=\int_0^{4\pi}I(\boldsymbol{x},\boldsymbol{r})d\omega.
\end{equation}
Here, the solid angle is indicated by \(\omega\).

\section{Scaling}
\label{scaling}
To construct non-dimensional bioconvection equations, scale the necessary factors, including lengths, pressure, fluid velocity, time, and cell concentration. Use standard parameters: \(H\), \(\mu \alpha_f/H^{2}\), \(\alpha_f/H\), \(H^{2}/\alpha_f\), and \(\bar{n}\). Non-dimensionalize temperature with \(\mathcal{T}=\frac{T^*-T_h^*}{\Delta T}\), where \(\alpha_f\) is the thermal diffusivity of water and \(\Delta T= T_l^*-T_h^*\). In the scenario of heating from above, $\Delta T$ is negative.

After non-dimensionalizing the governing equations (\ref{eqn:equation8})-(\ref{eqn:equation10}), become:
\begin{equation}
\label{eqn:equation13}
div(\boldsymbol{v})=0, 
\end{equation}
\begin{equation}
\label{eqn:equation14}
\frac{1}{P_r}\left(\frac{\partial}{\partial t}+\boldsymbol{v}\cdot\boldsymbol{\nabla}\right)\boldsymbol{v}=\nabla^2\boldsymbol{v}-\boldsymbol{\nabla} \mathcal{P}-n R_a \boldsymbol{k}-R_m \boldsymbol{k}+R_T \mathcal{T} \boldsymbol{k},
\end{equation}
\begin{equation}
\frac{\partial \mathcal{T} }{\partial t}+ \boldsymbol{v}\cdot \boldsymbol{\nabla} \mathcal{T}=\boldsymbol{\nabla}^2 \mathcal{T},
\end{equation}
\begin{equation}
\label{eqn:equation15}
\frac{\partial n}{\partial t}=-\boldsymbol{\nabla}\cdot\left(\frac{1}{Le}n U_s\bar{\boldsymbol{P}}+n\boldsymbol{v}-\frac{1}{Le}\boldsymbol{\nabla} n\right).
\end{equation}
Here, the basic-density Rayleigh number is $R_m=\frac{\varrho g H^3}{\mu \alpha_f}$, the bioconvection Rayleigh number is $R_a=\frac{\bar{n}\vartheta \Delta \rho g H^3}{\mu \alpha_f}$, the thermal Rayleigh number is $R_T=\frac{\beta \Delta T \varrho g H^3}{\mu \alpha_f}$, the Lewis number is $Le=\frac{\alpha_f}{D}$, the Prandtl number is $P_r=\frac{\mu}{\varrho \alpha_f}$, dimensionless swimming speed is $U_s=\frac{U_c H}{D}$, and the kinematic viscosity is $\nu=\mu/\rho$. The bioconvection Rayleigh number $R_a$ signifies the microorganisms' ability to induce bioconvection, whereas the thermal Rayleigh number $R_T$ quantifies the impact of the thermal factor on the entire system.

Non-dimensional boundary conditions for the stress-free top surface become:
\begin{align}
\label{eqn:equation17}
\boldsymbol{v}\cdot\boldsymbol{k}&=0, \quad \frac{\partial^2}{\partial x_3^2}(\boldsymbol{v}.\boldsymbol{k})=0, \quad \mathcal{T}=0, \quad \nonumber \\
&\left(\frac{1}{Le}n U_s\bar{\boldsymbol{P}}+n\boldsymbol{v}-\frac{1}{Le}\boldsymbol{\nabla} n\right)\cdot\boldsymbol{k}=0, \quad  \text{at} \quad x_3=1,
\end{align}
while for the rigid bottom surface
\begin{align}
\label{eqn:equation16}
\boldsymbol{v}&=0, \quad \boldsymbol{v}\times\boldsymbol{k}=0, \quad  \mathcal{T}=-1, \quad \nonumber\\
&\left(\frac{1}{Le}n U_s\bar{\boldsymbol{P}}+n\boldsymbol{v}-\frac{1}{Le}\boldsymbol{\nabla} n\right)\cdot\boldsymbol{k}=0, \quad \text{at} \quad x_3=0.
\end{align}

The non-dimensional radiative transfer equation becomes:
\begin{equation}
\frac{d I}{d r}+n \hbar  I(\boldsymbol{x},\boldsymbol{r})=0,
\end{equation}
where \(\hbar=\iota\bar{n}H\) is optical depth of the suspension. The non-dimensional boundary conditions for the top surface intensity are:
\begin{equation}
I(x_1,x_2,1,\alpha_1,\alpha_2)=I^0\delta(\boldsymbol{r}-\boldsymbol{r}^0) , \quad \pi/2\le\alpha_1\le\pi,
\end{equation}
and for the bottom surface intensity:
\begin{equation}
I(x_1,x_2,0,\alpha_1,\alpha_2)=0,\quad 0\le\alpha_1\le\pi/2.
\end{equation}

\section{Steady-state}
\label{sec3}
In the equilibrium state, we have $\boldsymbol{v}=0$, $\mathcal{P}=\mathcal{P}_b$, $n=n_b(x_3)$, $I=I_b(x_3,\alpha_1)$, $\mathcal{G}=\mathcal{G}_b(x_3)$, and $\mathcal{T}=\mathcal{T}_b$.

The steady-state radiative transfer equation becomes
\begin{equation}
\label{eqn:equation18}
\frac{\partial I_b}{\partial x_3}+\frac{\hbar n_b(x_3)}{\cos{\alpha_1}}I_b(x_3,\alpha_1)=0.
\end{equation}
Accompanied by the boundary condition
\begin{equation}
\label{eqn:equation19}
I_b(1,\alpha_1)=I^0\delta(\boldsymbol{r}-\boldsymbol{r}^0).
\end{equation}

Solving equation (\ref{eqn:equation18}), we obtain
\begin{equation}
\label{eqn:equation20}
I_b(x_3,\alpha_1)=C\exp{\left(\frac{-\hbar}{\cos{\alpha_1}}\int_1^{x_3} n_b(s)ds\right)}.
\end{equation}
Applying the boundary condition (\ref{eqn:equation19}), we get
\begin{equation}
\label{eqn:equation21}
I_b(x_3,\alpha_1)=I^0\delta(\boldsymbol{r}-\boldsymbol{r}^0)\exp{\left(\frac{-\hbar}{\cos{\alpha_1}}\int_1^{x_3} n_b(s)ds\right)}.
\end{equation}

The total intensity in the basic state $\mathcal{G}$ becomes
\begin{align}
\label{eqn:equation22}
\mathcal{G}_b(x_3)&=\int_0^{4\pi} I_b(x_3,\alpha_1)d\omega 
 =I^0\exp{\left(\hbar\int_1^{x_3} n_b(s)ds\right)}.
\end{align}

The cell conservation equation in the basic state transforms into
\begin{equation}
\label{eqn:equation24}
\frac{d n_b}{dx_3}-U_s T_b n_b=0,
\end{equation}
subject to the constraint
\begin{equation}
\label{eqn:equation25}
\int_0^1 n_b(x_3)dx_3=1.
\end{equation}

Expressing the radiative heat flux in the basic state as
\begin{align}
\label{eqn:equation26}
\boldsymbol{q}_b(x_3)&=\int_0^{4\pi}I_b(x_3,\alpha_1)\boldsymbol{r}d\omega =-I^0\exp{\left(\hbar\int_1^{x_3} n_b(s)ds\right)\boldsymbol{k}} \nonumber\
&=|\boldsymbol{q}_b|(-\boldsymbol{k}).
\end{align}

Consequently, the mean swimming orientation $\bar{\boldsymbol{P}}_b$ in the basic state is calculated by
\begin{equation}
\label{eqn:equation27}
\bar{\boldsymbol{P}}_b=-T(\mathcal{G}_b)\frac{\boldsymbol{q}_b}{|\boldsymbol{q}_b|}=T(\mathcal{G}_b)\boldsymbol{k}.
\end{equation}

The basic state temperature is given by
\begin{equation}
\mathcal{T}_b=x_3-1.
\end{equation}

\section{Linear stability analysis}
\label{sec4}

To investigate linear instability, we introduce an infinitesimal perturbation, denoted as $\epsilon (0 < \epsilon \ll 1)$, into the basic state. This perturbation is applied as follows:

\begin{align*}
    \boldsymbol{v} &= \boldsymbol{0} + \epsilon \boldsymbol{v}^*(x_1,x_2,x_3,t) + O(\epsilon^2),\\
    n &= n_b(x_3) + \epsilon n^*(x_1,x_2,x_3,t) + O(\epsilon^2),\\
    \mathcal{P} &= \mathcal{P}_b + \epsilon \mathcal{P}^* + O(\epsilon^2),\\
    \bar{\boldsymbol{P}} &= \bar{\boldsymbol{P}}_b + \epsilon\bar{\boldsymbol{P}}^* + O(\epsilon^2),\\
    \mathcal{G} &= \mathcal{G}_b + \epsilon\mathcal{G^*} + O(\epsilon^2),\\
    \mathcal{T} &= \mathcal{T}_b + \epsilon T' + O(\epsilon^2),
\end{align*}

where the perturbed fluid velocity is denoted as $\boldsymbol{v}^* = (u^*, v^*, w^*)$. Now, the linearized governing equations are expressed as follows:

\begin{equation}
\label{eqn:equation28a}
    \text{div}(\boldsymbol{v}^*) = 0,  
\end{equation}

\begin{equation}
\label{eqn:equation28b}
   \frac{1}{P_r}\frac{\partial \boldsymbol{v}^*}{\partial t} = \nabla^2\boldsymbol{v}^* - \boldsymbol{\nabla} \mathcal{P}^* - n^* R_a \boldsymbol{k} + R_T T' \boldsymbol{k},
\end{equation}

\begin{equation}
    \frac{\partial T'}{\partial t} - \boldsymbol{v}\cdot\boldsymbol{k} = \boldsymbol{\nabla}^2 T',
\end{equation}

\begin{equation}
\label{eqn:equation29}
    \frac{\partial n^*}{\partial t} = -\frac{dn_b}{dx_3}w^* + \frac{1}{Le}\nabla^2 n^* - \frac{1}{Le}U_s\boldsymbol{\nabla}\cdot(n^*\bar{\boldsymbol{P}}_b + n_b\bar{\boldsymbol{P}}^*).
\end{equation}

The total intensity, $\mathcal{G}$, is expressed as

\begin{align*}
     \mathcal{G} &= \mathcal{G}_b + \epsilon\mathcal{G^*} + O(\epsilon^2)\\
     &= I^0\exp{\left(\hbar\int_1^{x_3}(n_b(s)+\epsilon n^*+O(\epsilon^2))ds\right)}.
 \end{align*}

By accumulating $O(\epsilon)$ terms, the perturbed total intensity, $\mathcal{G}^*$, becomes

\begin{equation*}
    \mathcal{G}^* = I^0\left(\hbar\int_1^{x_3} n^* ds\right)\exp{\left(\hbar\int_1^{x_3} n_b(s)ds\right)}.
\end{equation*}

The mean swimming orientation is determined by
\begin{align}
\label{eqn:equation30}
    \bar{\boldsymbol{P}} &= \bar{\boldsymbol{P}}_b + \epsilon\bar{\boldsymbol{P}}^* + O(\epsilon^2) \nonumber \\
    &= T(\mathcal{G}_b + \epsilon \mathcal{G}^* + O(\epsilon^2))\boldsymbol{k}.
\end{align}

By accumulating $O(\epsilon)$ terms, the perturbed mean swimming orientation, $\bar{\boldsymbol{P}}^*$, is given by

\begin{equation}
\label{eqn:equation31}
    \bar{\boldsymbol{P}}^* = \mathcal{G}^*\frac{\partial T}{\partial \mathcal{G}_b}\boldsymbol{k}.
\end{equation}

Applying the curl operator twice to (\ref{eqn:equation28b}) and focusing on the vertical component yields 
\begin{align}
\label{eqn:equation34}
    \frac{1}{P_r}\frac{\partial}{\partial t}(\boldsymbol{\nabla}^2 w^*) &= \boldsymbol{\nabla}^4 w^* - R_a \boldsymbol{\nabla}^2 n^* + R_a\frac{\partial^2}{\partial x_3^2} n^* \nonumber \\
    &- R_T \boldsymbol{\nabla}^2 T' + R_T\frac{\partial^2}{\partial x_3^2} T'.
\end{align}
In equation (\ref{eqn:equation29}), we can express

\begin{eqnarray}
\label{eqn:equation32}
    \boldsymbol{\nabla}\cdot(n^*\bar{\boldsymbol{P}}_b + n_b\bar{\boldsymbol{P}}^*) &= 2\hbar n_b \mathcal{G}_b\frac{d T(\mathcal{G}_b)}{d \mathcal{G}_b} + T(\mathcal{G}_b)\frac{\partial n^*}{\partial x_3} + \hbar\frac{\partial}{\partial x_3}\left(n_b\mathcal{G}_b\frac{d T(\mathcal{G}_b)}{d\mathcal{G}_b}\right)\int_1^{x_3} n^* ds.
\end{eqnarray}

Therefore, equation (\ref{eqn:equation29}) can be rewritten as

\begin{align}
\label{eqn:equation33}
    &\frac{\partial n^*}{\partial t} - \frac{1}{Le}\boldsymbol{\nabla}^2 n^* - \frac{1}{Le}\hbar U_s\frac{\partial}{\partial x_3}\left(n_b\mathcal{G}_b\frac{d T(\mathcal{G}_b)}{d\mathcal{G}_b}\right)\int_{x_3}^1  n^* ds \nonumber\\
    &+ 2\frac{1}{Le}\hbar n_b \mathcal{G}_b\frac{d T(\mathcal{G}_b)}{d \mathcal{G}_b}n^*U_s + \frac{1}{Le}T(\mathcal{G}_b)\frac{\partial n^*}{\partial x_3}U_s = -w^*\frac{d n_b}{d x_3}.
\end{align}

Perturbed boundary conditions for the stress-free top surface become

\begin{align*}
    w^* &= 0, \quad \frac{\partial^2 w^*}{\partial {x_3}^2} = 0, \quad T' = 0, \quad U_s T_b n^* - \frac{\partial n^*}{\partial x_3} = 0, \quad \text{at} \quad x_3 = 1,
\end{align*}

while for the rigid bottom surface

\begin{align*}
    w^* &= 0, \quad \frac{\partial w^*}{\partial x_3} = 0, \quad T' = 0, \quad 
    \hbar U_s n_b \mathcal{G}_b \frac{d T_b}{d \mathcal{G}_b} \int_{x_3}^1 n^* d \bar{x}_3 \nonumber \\
    &- U_s T_b n^* + \frac{\partial n^*}{\partial x_3} = 0, \quad \text{at} \quad x_3 = 0.
\end{align*}

\section{Normal mode analysis}
 \label{normal}

Normal modes are decomposed from linearized governing equations via
\begin{equation*}
    w^* = \hat{w}(x_3)\exp{[\gamma t+i(a_1 x_1+a_2 x_2)]},
\end{equation*}
\begin{equation*}
    n^* = \hat{n}(x_3)\exp{[\gamma t+i(a_1 x_1+a_2 x_2)]},
\end{equation*}
\begin{equation*}
    T' = \Theta(x_3)\exp{[\gamma t+i(a_1 x_1+a_2 x_2)]}.
\end{equation*}
Here, $a_1$ and $a_2$ are wavenumbers in the $x_1$ and $x_2$ directions, and the resultant $a=\sqrt{a_1^2+a_2^2}$ is a horizontal wavenumber.

Thus, the governing equations (\ref{eqn:equation28a})-(\ref{eqn:equation29}) in normal modes become
\begin{align}
    \label{eqn:equation36}
    \frac{\gamma}{P_r}\left(\frac{d^2}{dx_3^2}-a^2\right) \hat{w}(x_3) -\left(\frac{d^2}{dx_3^2}-a^2\right)^2\hat{w}(x_3) 
    = a^2R_a\hat{n}(x_3)-a^2R_T\Theta(x_3),
\end{align}
\begin{equation}
    \label{eqn:equation37}
    \left(\gamma+a^2-\frac{d^2}{d x_3^2}\right)\Theta=\hat{w},
\end{equation}
\begin{eqnarray}
\label{eqn:equation38}
    &-\hbar U_s \frac{\partial}{\partial x_3}\left(n_b\mathcal{G}_b\frac{d T(\mathcal{G}_b)}{d\mathcal{G}_b}\right)\int_{x_3}^1  \hat{n}d \bar{x}_3+( \gamma Le+a^2) \hat{n} \nonumber \\
    &+2\hbar n_b \mathcal{G}_b\frac{d T_b}{d \mathcal{G}_b}U_s \hat{n}+U_s T_b \frac{d \hat{n}}{d x_3}
    -\frac{d^2 \hat{n}}{d x_3^2}
    =-Le \frac{d n_b}{d x_3}\hat{w}.
\end{eqnarray}

The associated boundary conditions in the normal modes for the stress-free top surface become
\begin{align*}
    \hat{w}(x_3)=0, \quad \frac{d^2\hat{w}(x_3)}{d x_3^2}=0, \quad \Theta=0, \quad 
    U_s T_b \hat{n}-\frac{d \hat{n}}{d x_3}=0  \quad \text{at} \quad x_3=1,
\end{align*}
while for the rigid bottom surface
\begin{align*}
    \hat{w}(x_3)=0, \quad \frac{\hbar\hat{w}(x_3)}{d x_3}=0, \quad \Theta=0, \quad
    \hbar U_s n_b \mathcal{G}_b \frac{d T_b}{d \mathcal{G}_b} \int_{x_3}^1 \hat{n} d \bar{x}_3-U_s T_b \hat{n}+\frac{d \hat{n}}{d x_3}=0
    \quad \text{at} \quad x_3=0.
\end{align*}

Define a new variable
\begin{equation}
\label{eqn:equation39}
    N(x_3)=\int_{x_3}^1 \hat{n} d \bar{x}_3,
\end{equation}
so that the system of equations becomes
\begin{align}
\label{eqn:equation40}
    \frac{d^4 \hat{w}}{d x_3}-\left(2 a^2+ \frac{\gamma}{P_r}\right)\frac{d^2 \hat{w}}{d x_3^2}+a^2\left(a^2+\frac{\gamma}{P_r}\right)\hat{w}
    =a^2 R_a \frac{d N}{d x_3}+a^2 R_T \Theta,
\end{align}
\begin{equation}
\label{eqn:equation41}
    \left(\gamma+a^2-\frac{d^2}{d x_3^2}\right)\Theta=\hat{w},
\end{equation}
\begin{align}
\label{eqn:equation42}
    &\frac{d^3 N}{d x_3^3}-U_s T_b \frac{d^2 N}{d x_3^2}-\left(\gamma Le+a^2+2\hbar U_s n_b \mathcal{G}_b\frac{d T_b}{d \mathcal{G}_b}\right)\frac{d N}{d x_3} \nonumber \\
    &- \hbar U_s \frac{d}{d x_3}\left(n_b \mathcal{G}_b\frac{d T_b}{d \mathcal{G}_b}\right)N=-Le\frac{d n_b}{d x_3}\hat{w}.
\end{align}

Also, boundary conditions for the stress-free top surface become
\begin{align}
    \label{eqn:equation45}
    \hat{w}(x_3)=0, \quad \frac{d^2\hat{w}(x_3)}{d z^2}=0, \quad \Theta=0, \quad 
    U_s T_b \frac{d N}{d x_3}-\frac{d^2 N}{d x_3^2}=0 \quad \text{at} \quad x_3=1,
\end{align}
while for the rigid bottom surface
\begin{align}
    \label{eqn:equation43}
    \hat{w}(x_3)=0, \quad \frac{d\hat{w}(x_3)}{d x_3}=0, \quad \Theta=0, \quad
    \hbar U_s n_b \mathcal{G}_b \frac{d T_b}{d \mathcal{G}_b} N+U_s T_b \frac{d N}{d x_3}
    -\frac{d^2 N}{d x_3^2}=0 \quad \text{at} \quad x_3=0,
\end{align}

and 
\begin{equation}
\label{eqn:equation46}
    N(x_3)=0 \quad \text{at} \quad x_3=1.
\end{equation}

A set of governing equations and boundary conditions (\ref{eqn:equation40})-(\ref{eqn:equation46}) constitute a system of phototactic thermal-bioconvection equations. Without thermal convection, the governing equations return to the form presented by Vincent and Hill \cite{ref9}.

\begin{figure*}
    \centering
    \includegraphics[width=15cm, height=15cm ]{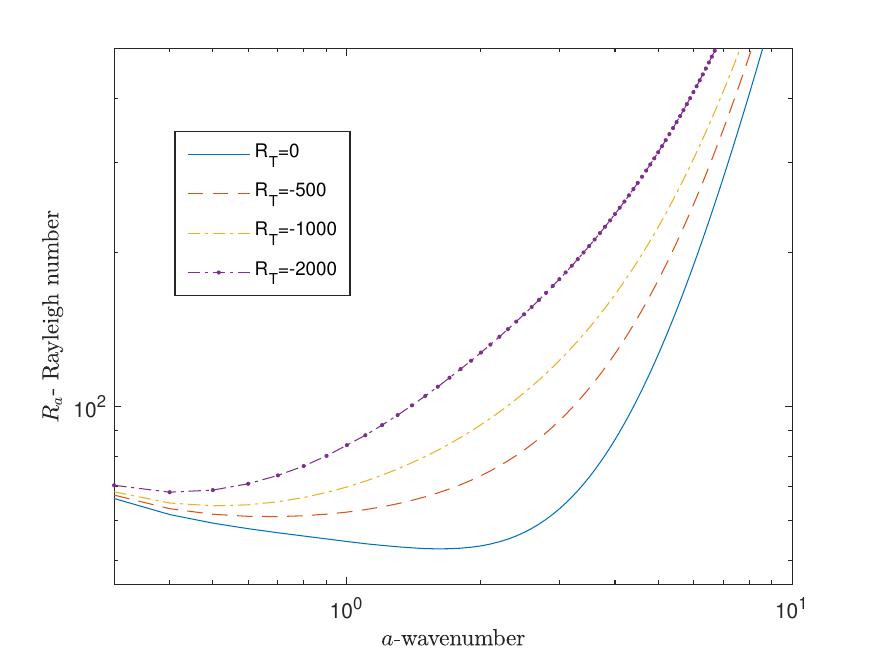}
    \caption{Neutral curves fixed parameters $\hbar=0.5$, $U_s=10$, $\mathcal{G}_c=0.68$, and $Le=4$ as $R_T$ is varied.}
   \label{fig:10v0.5k0.68Ic.pdf}
 \end{figure*}
 \begin{figure*}
    \centering
    \includegraphics[width=15cm, height=15cm ]{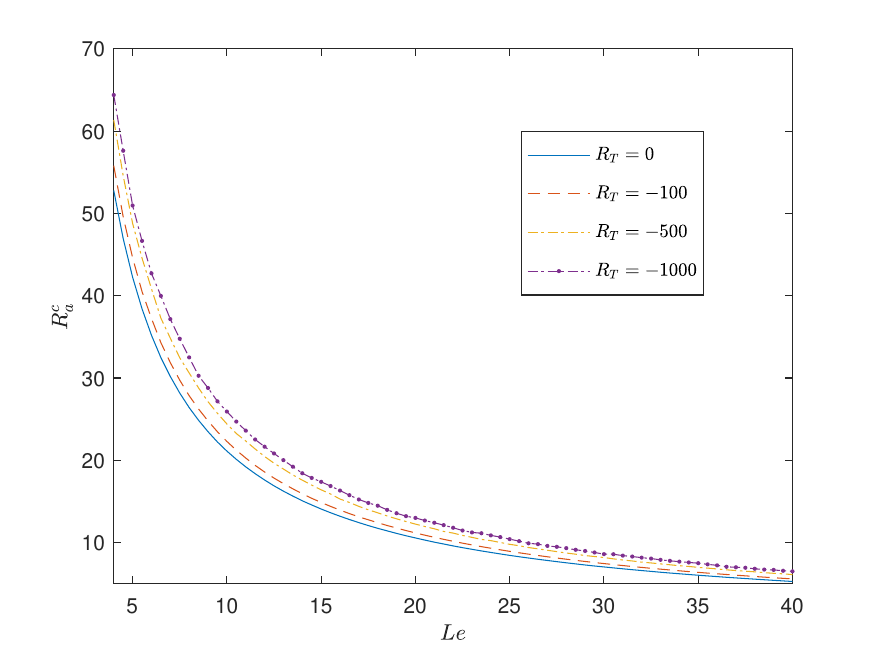}
    \caption{The relationship between $R_a^c$ and $Le$ for $\hbar=0.5$, $U_s=10$, and $\mathcal{G}_c=0.68$ as $R_T$ is varied.  }
   \label{fig:effofLe.pdf}
 \end{figure*}
\begin{figure*}
    \centering
    \includegraphics[width=15cm, height=15cm ]{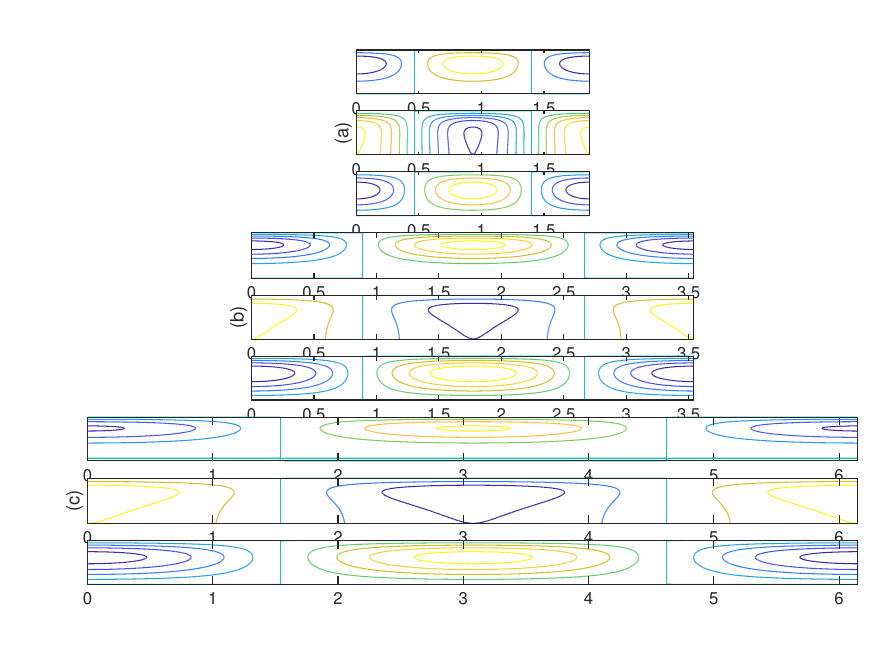}
    \caption{Streamlines (upper), isoconcentration (middle) and isotherm (lower) for fixed parameters $\hbar=0.5$, $U_s=10$, $Le=4$, $R_a=70$ and $\mathcal{G}_c=0.68$; (a) $R_T=0$, (b) $R_T=-500$, (c) $R_T=-1000$.  }
   \label{fig:stream.pdf}
 \end{figure*}
 \begin{figure*}
    \centering
    \includegraphics[width=15cm, height=15cm ]{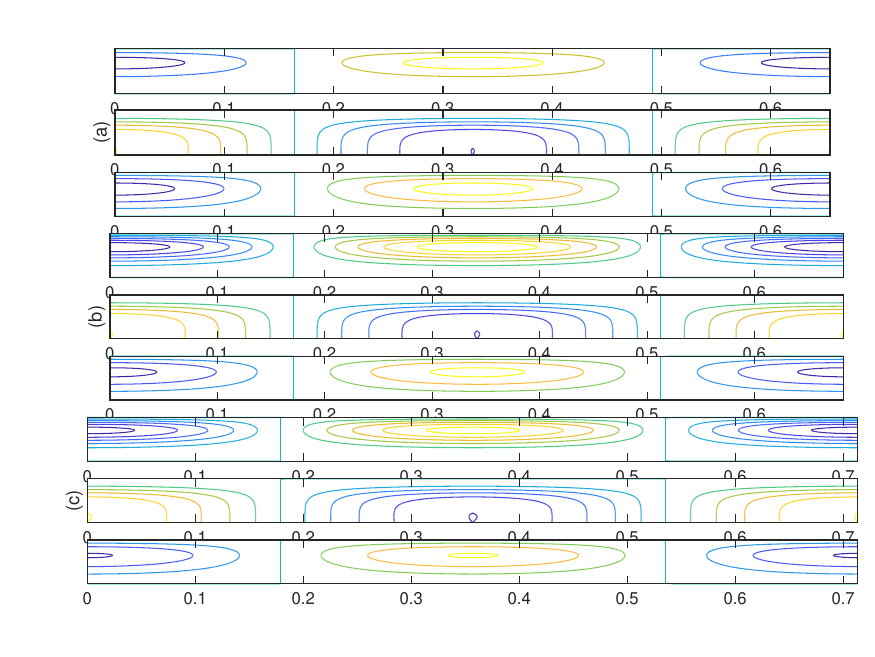}
    \caption{Streamlines (upper), isoconcentration (middle) and isotherm (lower) for fixed parameters $\hbar=0.5$, $U_s=10$, $Le=40$, $R_a=70$ and $\mathcal{G}_c=0.68$; (a) $R_T=0$, (b) $R_T=-500$, (c) $R_T=-1000$.  }
   \label{fig:stream40Le.pdf}
 \end{figure*}
 \begin{figure*}
    \centering
    \includegraphics[width=16cm, height=11cm ]{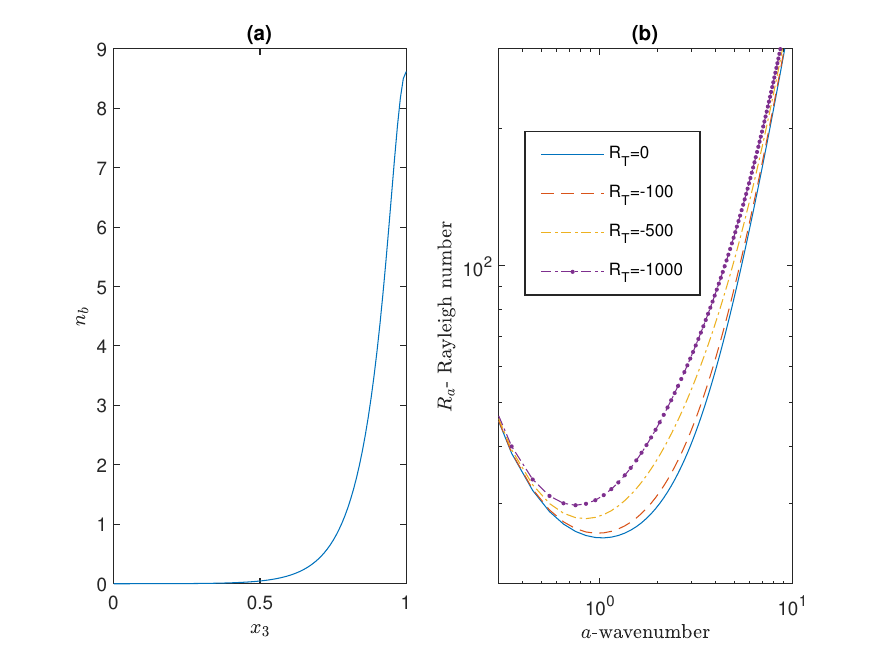}
    \caption{(a) The basic state concentration and (b) neutral curves fixed parameters $\hbar=0.5$, $U_s=15$, $\mathcal{G}_c=0.8$, and $Le=4$ as $R_T$ is varied. }
   \label{fig:15v0.5k0.8Ic.pdf}
 \end{figure*}
\begin{figure*}
    \centering
    \includegraphics[width=16cm, height=11cm ]{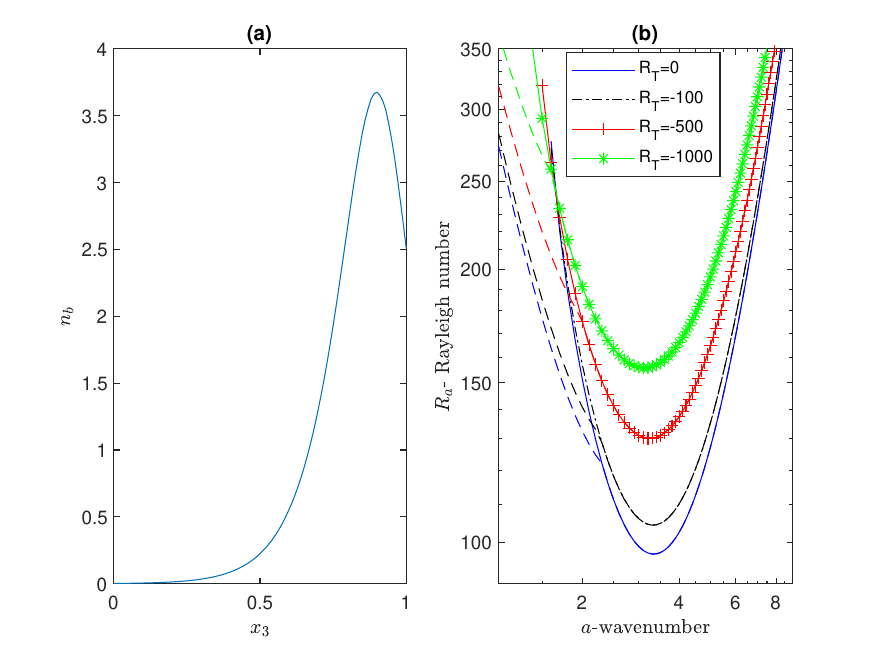}
    \caption{(a) The basic state concentration and (b) neutral curves fixed parameters $\hbar=0.5$, $U_s=15$, $\mathcal{G}_c=0.68$, and $Le=4$ as $R_T$ is varied. The dotted lines represented the oscillatory branch.}
   \label{fig:15v0.5k0.68Ic.pdf}
 \end{figure*}
\begin{figure*}
    \centering
    \includegraphics[width=16cm, height=11cm ]{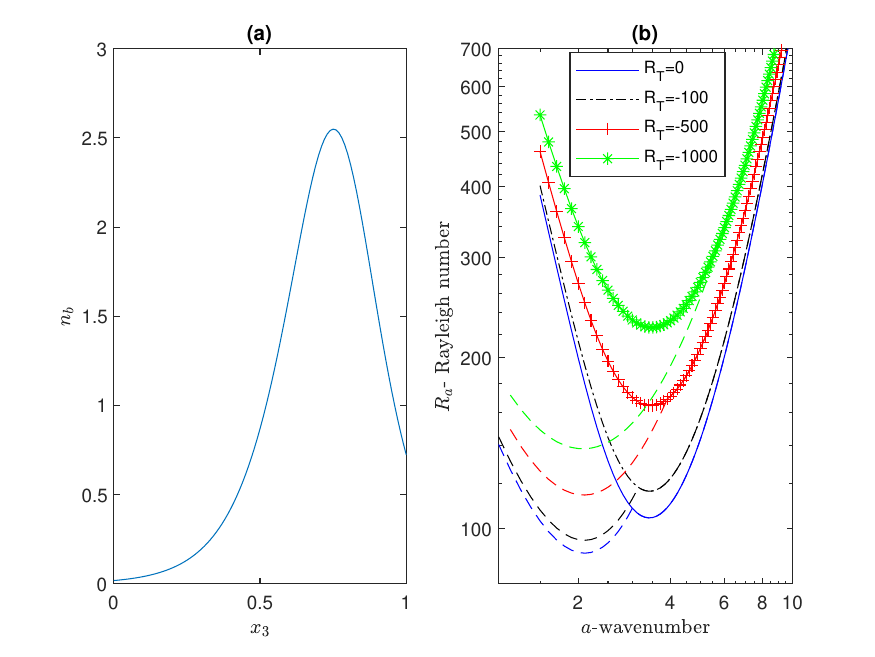}
    \caption{(a) The basic state concentration and (b) neutral curves fixed parameters $\hbar=0.5$, $U_s=15$, $\mathcal{G}_c=0.65$, and $Le=4$ as $R_T$ is varied.  The dotted lines represented the oscillatory branch.}
   \label{fig:15v0.5k0.65Ic.pdf}
 \end{figure*}
 \begin{figure*}
    \centering
    \includegraphics[width=16cm, height=11cm ]{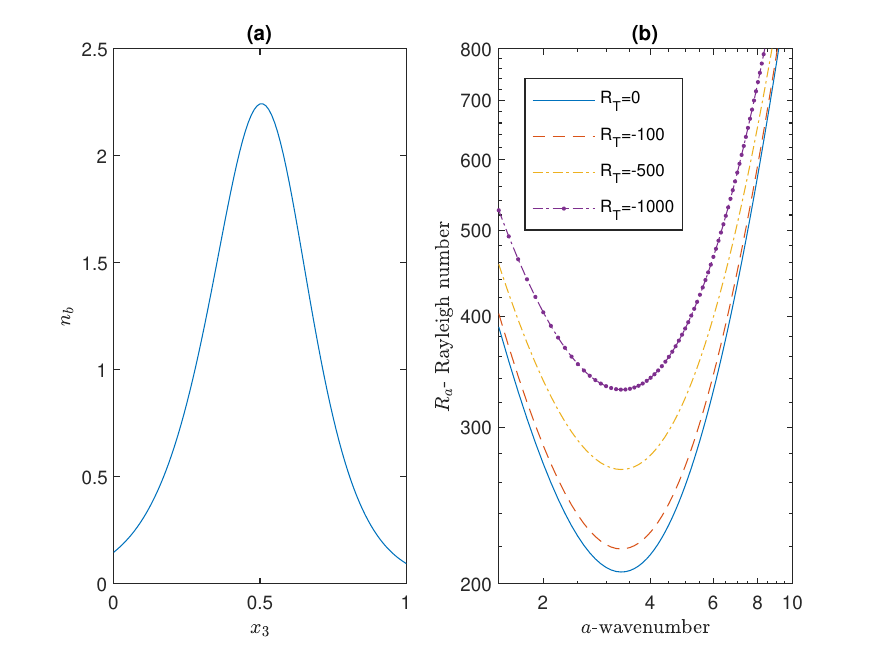}
    \caption{(a) The basic state concentration and (b) neutral curves fixed parameters $\hbar=0.5$, $U_s=15$, $\mathcal{G}_c=0.63$, and $Le=4$ as $R_T$ is varied. }
   \label{fig:15v0.5k0.63Ic.pdf}
 \end{figure*}
\begin{table*}
\caption{Common suspension parameters for the phototactic microorganism \textit{Chlamydomonas}\cite{ref9,ref13,ref17,kumar2023effect,zhao2018linear}.}
\begin{tabular}{ p{8cm} p{5cm}}
\hline
\hline
Scaled average swimming speed & $U_s=20H$\\
Kinematic viscosity &$\nu=10^{-2}$cm$^2$/s\\
Prandtl number &$P_r=5$ \\
Average concentration &$\bar{n}=10^6$cm$^{-3}$\\
Average cell swimming speed &$U_c =10^{-2}$cm/s\\
Cell volume &$\vartheta=5\times10^{-10}$cm$^3$\\
Cell diffusivity &$D=5\times10^{-5}-5\times10^{-4}$cm$^2$/s\\
Thermal diffusivity&$\alpha_f=2\times 10^{-3}$cm$^2$/s\\
Volumetric thermal expansion coefficient&$\beta=3.4\times 10^{-3}$K$^{-1}$\\
Temperature of upper wall&$T_h^*=300$ K\\
Temperature difference& $\Delta T=-1$K\\
Ratio of cell density &$\Delta \varrho/\varrho=5\times 10^{-2}$ \\
Cell radius &$10^{-3}$cm\\
\hline
\hline
\end{tabular}\\
\label{tab:table1}
\end{table*}

 \begin{figure*}
    \centering
    \includegraphics[width=15cm, height=15cm ]{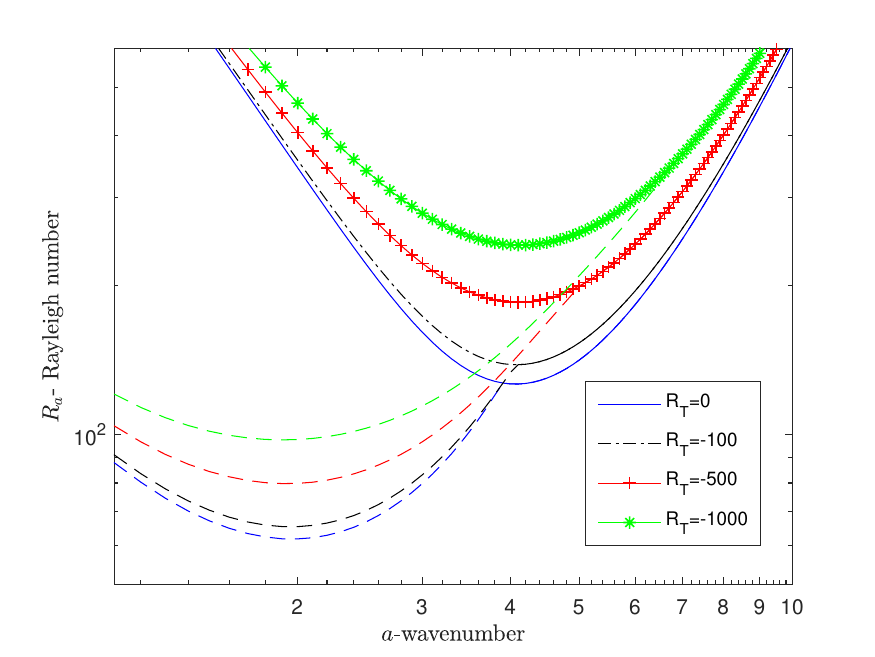}
    \caption{Neutral curves for $\hbar=1.0$, $U_s=15$, $\chi=-0.485$, $\mathcal{G}_c=0.51$, $Le=4$ as $R_T$ is varied. The dotted lines represented the oscillatory branch.}
   \label{fig:15v1k.pdf}
 \end{figure*}

\begin{figure*}
    \centering
    \includegraphics[width=16cm, height=20cm ]{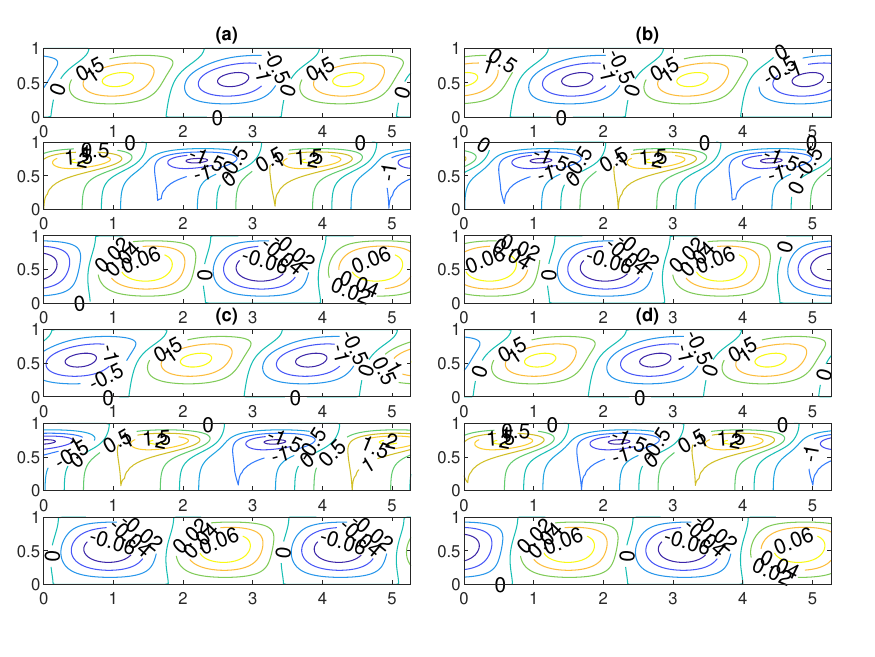}
    \caption{Flow pattern produced by the perturbed velocity (upper), concentration (middle) and temperature (lower) throughout one cycle of oscillation for  $\hbar=1.0$, $U_s=15$, $\chi=-0.485$, $\mathcal{G}_c=0.51$, $Le=4$, $R_T=-500$, $R_a^c=79.78$, $a_c=1.9$, $Im(\gamma)=12.98$, (a) $t=0$, (b) $t=0.16$, (c) $t=0.32$, and (d) $t=0.48$.}
   \label{fig:periodic.pdf}
 \end{figure*}
 
\begin{figure*}
    \centering
    \includegraphics[width=16cm, height=11cm ]{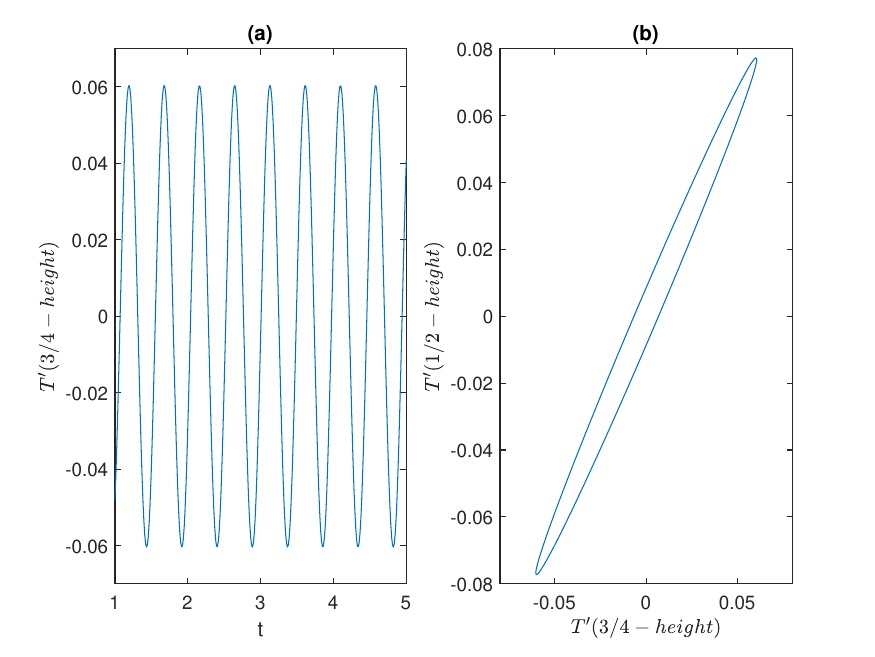}
    \caption{(a) Time-evolving perturbed temperature $T'$ and (b) phase diagram for $\hbar=1.0$, $U_s=15$, $\chi=-0.485$, $\mathcal{G}_c=0.51$, $Le=4$, $R_T=-500$.}
   \label{fig:timeinv+limitcycleRt.pdf}
 \end{figure*}
\begin{figure*}
    \centering
    \includegraphics[width=16cm, height=11cm ]{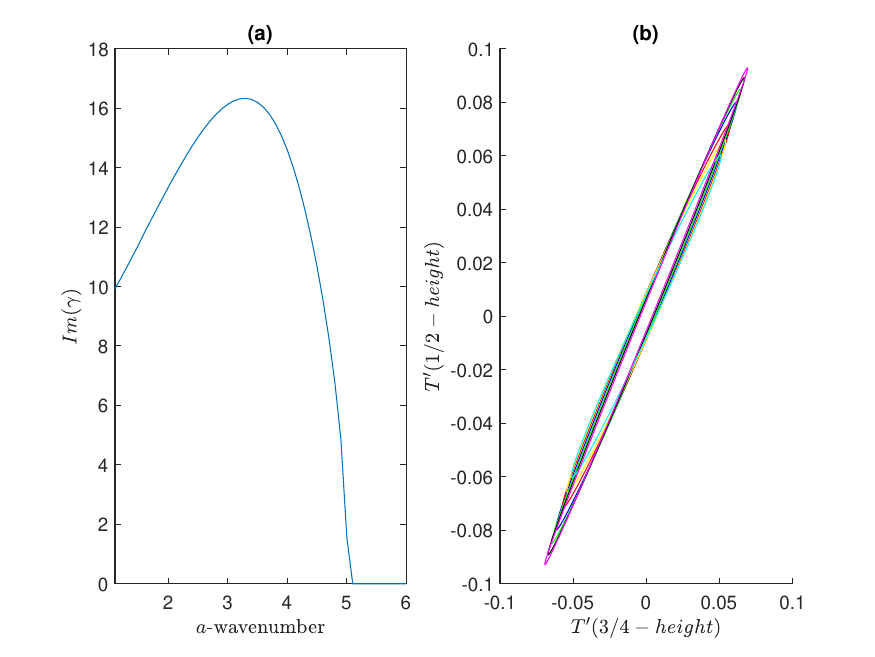}
    \caption{(a) The positive frequency as a function of the wavenumber $a$ corresponding to the oscillatory branch of the neutral curve for the thermal Rayleigh number $R_t=-500$ and (b) the corresponding phase portrait of the perturbed temperature $T'$. Fixed parameter values are $\hbar=1.0$, $U_s=15$, $\chi=-0.485$, $\mathcal{G}_c=0.51$, $Le=4$.}
   \label{fig:limitcycleRtImtest.pdf}
 \end{figure*}

 \section{Numerical results}
\label{sec5}
The system of ordinary differential equations (ODEs) given by equations (\ref{eqn:equation40}) to (\ref{eqn:equation42}) forms a ninth-order system accompanied by nine boundary conditions (\ref{eqn:equation45}) to (\ref{eqn:equation46}). To address coupled ODEs, the shooting technique, a numerical method, is employed as part of the analysis approach. The bvp4c solver in MATLAB is utilized for computation \cite{shampine2003solving}. It is noteworthy that when the thermal Rayleigh number $R_T$ is altered, the fluid medium exhibits characteristics akin to the pure bioconvection problem. Similarly, when $R_a$ is modified, it corresponds to the manifestation of B$\acute{e}$nard thermoconvection. This work primarily focuses on investigating the linear stability of the basic state, emphasizing the illustration of neutral curves. These curves comprise points where the real component of the growth rate $Re(\gamma)$ is zero. These zero points are pivotal in stability analysis. The existence of an overstable or oscillatory solution is possible if the imaginary component of the growth rate $Im(\gamma)$ along the neutral curve is nonzero. Conversely, according to the concept of exchange of stabilities, perturbations to the basic state are considered stationary if $Im(\gamma)$ is zero along this curve \cite{ref12}. In the time dependence expressed as $\exp(\gamma t) = \exp[(Re(\gamma)+i Im(\gamma))t]$, both the growth rate $Re(\gamma)$ and a hypothetical oscillation described by $Im(\gamma)/2\pi$ are included. For stable layering, the wave vectors have a negative value for the $Re(\gamma)$ parameter. However, just above the convective threshold, a narrow range comes into play where $Re(\gamma)$ becomes positive. The branch of the solution that exhibits the most significant instability is the one where the bioconvection Rayleigh number $R_a$ takes on its lowest possible value, denoted as $R_a^c$. This critical solution, labeled as $(a^c, R_a^c)$, is recognized as the most unstable solution. To align our model with prior work on phototactic bioconvection, we consider phototactic microorganisms, analogous to \textit{Chlamydomonas}. The necessary parameters for this investigation are determined by consulting previous studies \cite{ref9, ref13, ref17, kumar2023effect} (refer to Table \ref{tab:table1}). The radiation parameters are consistent with those in Ref. \cite{ref14}. The optical depth can vary from $0.25$ to $1$ for a suspension depth of $0.5$ cm, and $U_s = 10$ is the scaled swimming speed corresponding to this range. Similarly, the value of $U_s$ for a suspension with a depth of $1.0$ cm is $20$. To maintain compatibility with other phototactic bioconvection models, $I^0 = 0.8$ is kept constant throughout. Previous research has shown that the critical intensity $\mathcal{G}_c$ is linked to the parameter $\chi$ in a way that $-1.1 \leq \chi \leq 1.1$ leads to $0.3 \leq \mathcal{G}_c \leq 0.8$ \cite{ref17, kumar2023effect}.

Stability diagrams for fixed parameters $\hbar=0.5$, $U_s=10$, $\mathcal{G}_c=0.68$ and $Le=4$, are shown in Figure \ref{fig:10v0.5k0.68Ic.pdf}, where the bioconvection Rayleigh number $R_a$ is plotted against wavenumber $a$ and as a function of the thermal Rayleigh number $R_T$. We consider the case of thermal gradients heated from above (or cooled from below) with a negative thermal Rayleigh number. The Rayleigh number $R_a$, influenced by biological parameters, varies as a function of the wave number $a$, reaching its minimum at a critical wave number $a^c$. Decreasing the thermal Rayleigh number $R_T$ results in an elevation of the critical bioconvection Rayleigh number $R_a^c$, enhancing the stability of the system. When the vector of the temperature gradient aligns with the gravity vector and the thermal Rayleigh number $R_T$ is negative, the Rayleigh-B$\acute{e}$nard convection theory anticipates the absence of thermal convection due to heating from above.  Consequently, the predominant mechanism governing the system is microorganism motion convection, commonly known as bioconvection. It is noteworthy that temperature continues to play a significant role in stabilizing bioconvection-induced convection through the parameter critical bioconvection Rayleigh number $R_a$. The connection between the critical Rayleigh number $R_a^c$ and Lewis number $Le$, with parameters $\hbar=0.5$, $U_s=10$, and $\mathcal{G}_c=0.68$, is illustrated in Figure \ref{fig:effofLe.pdf} while varying the thermal Rayleigh number $R_T$. We observe that increasing the Lewis number leads to a decrease in the critical bioconvection Rayleigh number. This implies that the system becomes more unstable with a higher Lewis number.

Figures \ref{fig:stream.pdf}-\ref{fig:stream40Le.pdf} depict the streamlines, isoconcentration, and isotherm patterns, with fixed parameters Rt at values of $\hbar=0.5$, $U_s=10$, $Le=4$, $R_a=70$ and $\mathcal{G}_c=0.68$; (a) $R_T=0$, (b) $R_T=-500$, (c) $R_T=-1000$. The graphics demonstrate how the thermal Rayleigh number $R_T$ significantly affects both the distribution of microorganisms in the system and the patterns of convection cells. In Figure \ref{fig:stream.pdf}, we observe that the pattern wavelength $\lambda$ is $1.86$ for $R_T=0$. However, as $R_T$ is decreased to $-1000$, the pattern wavelength increases to $6.15$, indicating a substantial influence of the thermal Rayleigh number $R_T$ reduction on the pattern wavelength. This change suggests a notable impact on the system. Interestingly, when the Lewis number $Le$ is increased up to $40$, increments in $R_T$ do not exhibit significant variations in the pattern wavelength.

Figures \ref{fig:15v0.5k0.8Ic.pdf}-\ref{fig:15v0.5k0.63Ic.pdf} demonstrate the impact of the thermal Rayleigh number $R_T$ on both the fundamental concentration and neutral curves, considering various values of the critical total intensity $\mathcal{G}_c$. Here the fixed parameters are $\hbar=0.5$, $U_s=15$, and $Le=4$. A horizontally concentrated sublayer is produced due to the highest basic concentration occurring near the top of the suspension when $G_c=0.8$ (see Figure \ref{fig:15v0.5k0.8Ic.pdf}(a)). The region above the sublayer is gravitationally stable, while the region below it is unstable. In this case, the highest basic state concentration, $n_{b_{(max)}}=8.61$ is found at $x_3=1$. Under these conditions, when the thermal Rayleigh number $R_T$ is varied from $0$ to $-1000$, the critical bioconvection Rayleigh number increases (see Figure \ref{fig:15v0.5k0.8Ic.pdf}(b)). Additionally, in this scenario, only the stationary branch is identified. There is no oscillatory branch found. As the critical total intensity $\mathcal{G}_c$ decreases to $0.68$, the position of the highest basic state concentration, or the sublayer, moves away from the top of the suspension. Simultaneously, the maximum basic concentration value also decreases (see Figure \ref{fig:15v0.5k0.68Ic.pdf}(a)). In this case, the highest basic state concentration, $n_{b_{(max)}}=3.37$ is found at $x_3=0.9$. For $R_T=0$, an oscillatory branch converges with the stationary branch at approximately $a=2.4$. The oscillatory branch is present for $a\le2.4$, but the most unstable mode is situated on the stationary branch (see Figure \ref{fig:15v0.5k0.68Ic.pdf}(b)). As the thermal Rayleigh number $R_T$ decreases to $-1000$, the critical Rayleigh number $R_a^c$ increases and the suspension becomes more stable.

When $\mathcal{G}_c$ reaches to $0.65$, the highest basic state concentration, $n_{b_{(max)}}=2.54$ is found at $x_3=0.76$ (see Figure \ref{fig:15v0.5k0.65Ic.pdf}(a)). When the thermal Rayleigh number $R_T$ is set to $0$, a single oscillatory branch emerges from the stationary branch of the neutral curve, occurring approximately at $a=3.1$. In this state, the most unstable solution is identified at $(R_a^c, a^c) = (90.54, 2.1)$ on the oscillatory branch, and the perturbation to the basic state becomes overstable. The positive frequency, $Im(\gamma)$, associated with the most unstable solution has an approximate value of $6.81$, resulting in a period of overstability that is approximately $0.92$ units. On the thermal Rayleigh number $R_T$ varying, the most unstable solution remains on the overstable branch (see Figure \ref{fig:15v0.5k0.65Ic.pdf}(b)). Oscillatory instability typically arises when the instability results from the interplay of two or more competing processes. One of these processes pushes the system away from its basic state, while the other attempts to restore the system to its basic state \cite{ref12}. In a phototactic suspension, three distinct processes come into play. The gravitationally stable region situated above the layer of maximum concentration acts as an inhibitor of convection, whereas the region below supports it. Phototaxis plays a dual role; it can either inhibit or support convection. The presence of these competing processes results in the observation of oscillatory solutions. When $\mathcal{G}_c$ reaches $0.63$, the highest basic state concentration, $n_{b_{(max)}}=2.24$ is found at $x_3=0.52$ (see Figure \ref{fig:15v0.5k0.63Ic.pdf}(a)). In this scenario, the sublayer is positioned approximately around the mid-height of the suspension. Moreover, the region supporting positive phototaxis, which supports convection, decreases. As a result, the suspension becomes more stable for low values of critical total intensity. In this situation, the oscillatory branch of the neutral curve disappears entirely, and the most unstable solution persists within the stationary branch of the neutral curve, resulting in a transition to a stationary solution (see Figure \ref{fig:15v0.5k0.63Ic.pdf}(b)).

Neutral curves for the given parameters $\hbar=1.0$, $U_s=15$, $\chi=-0.485$, $\mathcal{G}_c=0.51$, and $Le=4$ are illustrated as the thermal Rayleigh number $R_T$ is varied. When the thermal Rayleigh number $R_T$ is set to zero, a solitary oscillatory branch emerges through a bifurcation from the stationary branch of the corresponding neutral curve, occurring approximately at $a=3.9$. The most unstable mode on the oscillatory branch of the neutral curve, and consequently, the bioconvective solution, is characterized by being overstable. With a thermal Rayleigh number $R_T$ set to $-500$, a lone oscillatory branch undergoes bifurcation from the corresponding stationary branch of the neutral curve at approximately $a=5.1$. However, the oscillatory branch maintains the most unstable bioconvective solution. As a result, the onset of overstability occurs at $a^c=1.9$ and $R_a^c=79.78$. Two complex conjugate eigenvalues, denoted as $\gamma=0\pm12.98$, are identified at this juncture. The observed transition is recognized as a Hopf bifurcation. The bioconvective flow patterns associated with the complex conjugate pair of eigenvalues exhibit mirror images of each other. The period of oscillation is measured as $2\pi/Im(\gamma)=0.48$ units. The bioconvective fluid motions attain full nonlinearity within a timescale significantly shorter than the anticipated period of overstability. Therefore, the $w^*$, $n^*$, and $T'$ perturbed eigenmodes can be used to see the convection cells and flow patterns during one oscillation cycle (see Figure \ref{fig:periodic.pdf}). This observation indicates that a traveling wave solution is propagating towards the left side of the figure. Figure \ref{fig:timeinv+limitcycleRt.pdf} illustrates the predicted time-evolving perturbed temperature component $T'$ in (a) and its corresponding phase diagram in (b) at critical wavenumber $a^c=1.9$. It is noteworthy that the period of oscillation, denoted as $2\pi/Im(\gamma)$, serves as the bifurcation/control parameter. Consequently, the destabilization of bioconvective flow gives rise to a limit cycle, represented as an isolated cycle in the phase diagram (see Figure \ref{fig:timeinv+limitcycleRt.pdf} (b)). The emergence of a limit cycle due to flow destabilization is identified as the Hopf bifurcation through bifurcation analysis. According to linear stability theory, the supercritical nature of this Hopf bifurcation ultimately leads it into a stable limit cycle. In Figure \ref{fig:limitcycleRtImtest.pdf}(a), the dependence of the positive frequency on the wavenumber $a$ is presented for the oscillatory branch of the neutral curve when $R_T=-500$. Despite the eigenvalues $\gamma$ appearing in complex conjugate pairs, only the positive component of the frequency, $Im(\gamma)$, is depicted in Figure \ref{fig:limitcycleRtImtest.pdf}(a). As the frequency approaches zero, the oscillatory mode of the disturbance transitions into a stationary state at the onset of bioconvection. The period of oscillation is denoted as $2\pi/Im(\gamma)$, where $Im\ne0$ serves as the control (bifurcation) parameter, aiding in the visualization of phase portraits and bifurcation diagrams. This parameter can be modified by multiplying it with an integer ($m2\pi/Im(\gamma)$, where $m$ is an integer), leading to periodic orbits or limit cycles in the behavior of the bioconvective system.

Figure \ref{fig:limitcycleRtImtest.pdf}(b) illustrates that the bifurcation diagram displays a degrading orbit/non-isolated cycle with sequentially diminishing radii for each orbit. This suggests that the bioconvective flow regime exhibits damped oscillations as long as the corresponding frequency has a non-zero value, specifically when $a<5.1$ (see Figure \ref{fig:limitcycleRtImtest.pdf}(a)). When the frequency tends to zero ($a^5.1$), the bioconvective flow regime undergoes a transition to steady convective motion.
Even when reducing the thermal Rayleigh number $R_T$ to $-1000$, the most unstable mode continues to persist on the oscillatory branch.

\section{Conclusion}
\label{sec6}

This study highlights the significant impact of upper-surface heating and collimated irradiation on phototactic bioconvection in non-scattering suspensions. The linear stability of the system was explored, focusing on neutral curves where the real growth rate is zero. The existence of overstable or oscillatory solutions was contingent on the imaginary growth rate along these curves. Variations in the thermal Rayleigh number, the Lewis number, and the critical total intensity were explored. Lower thermal Rayleigh numbers led to increased stability, and higher Lewis numbers and the critical total intensity resulted in a decreased critical bioconvection Rayleigh number, indicating heightened instability. Negative values of the thermal Rayleigh number indicate a scenario where the temperature difference is negative, meaning the upper plate is at a higher temperature than the lower plate. In the absence of bioconvection, the layer with a negative temperature gradient would remain stable. However, the upswimming motion of microorganisms introduces two conflicting effects in this situation: the destabilizing influence of microorganism upswimming on density stratification and the stabilizing impact of the temperature gradient on density stratification. The visualizations presented streamlines, isoconcentration, and isotherm patterns, which revealed their sensitivity to changes in thermal Rayleigh numbers and limited variations in pattern wavelength for greater Lewis numbers. Oscillatory branches, bifurcations, and Hopf bifurcations were identified in neutral curves, underscoring the significance of control parameters. The study enhances our understanding of bioconvection dynamics, particularly in phototactic microorganism systems. The results lay the foundation for further exploration and applications in related fields.

To make sure that the proposed phototaxis model is correct, it is important to check the theoretical predictions against measured experimental observations of bioconvection in a solution of phototactic algae. Unfortunately, as of now, there are no available statistics for such comparisons. The challenge lies in finding a microbe species that is predominantly phototactic, as most algae in natural environments also exhibit gravitactic or gyrotactic behavior \cite{ref8}. To gain further insights into these phenomena, additional experimentation is required to determine optical depths, phototaxis functions, and diffusion coefficients. These experiments will help unravel the complexities of bioconvection in phototactic algae solutions. We can learn more about how temperature affects suspensions through different mechanisms with the help of the proposed model. These mechanisms include isotropic scattering, anisotropic scattering, diffused irradiation, and oblique irradiation.

\section{Data availability}
The data that supports the findings of this study is in the article. 
\section{Declaration of Interests}
The authors state that they have no conflicting interests. 
\section{Acknowledgements}
The University Grants Commission financially supports this study, Grants number 191620003662, New Delhi (India).

\nocite{*}
\bibliography{aipsamp}

\providecommand{\noopsort}[1]{}\providecommand{\singleletter}[1]{#1}%
\begin{thebibliography}{29}%
\makeatletter
\providecommand \@ifxundefined [1]{%
 \@ifx{#1\undefined}
}%
\providecommand \@ifnum [1]{%
 \ifnum #1\expandafter \@firstoftwo
 \else \expandafter \@secondoftwo
 \fi
}%
\providecommand \@ifx [1]{%
 \ifx #1\expandafter \@firstoftwo
 \else \expandafter \@secondoftwo
 \fi
}%
\providecommand \natexlab [1]{#1}%
\providecommand \enquote  [1]{``#1''}%
\providecommand \bibnamefont  [1]{#1}%
\providecommand \bibfnamefont [1]{#1}%
\providecommand \citenamefont [1]{#1}%
\providecommand \href@noop [0]{\@secondoftwo}%
\providecommand \href [0]{\begingroup \@sanitize@url \@href}%
\providecommand \@href[1]{\@@startlink{#1}\@@href}%
\providecommand \@@href[1]{\endgroup#1\@@endlink}%
\providecommand \@sanitize@url [0]{\catcode `\\12\catcode `\$12\catcode `\&12\catcode `\#12\catcode `\^12\catcode `\_12\catcode `\%12\relax}%
\providecommand \@@startlink[1]{}%
\providecommand \@@endlink[0]{}%
\providecommand \url  [0]{\begingroup\@sanitize@url \@url }%
\providecommand \@url [1]{\endgroup\@href {#1}{\urlprefix }}%
\providecommand \urlprefix  [0]{URL }%
\providecommand \Eprint [0]{\href }%
\providecommand \doibase [0]{http://dx.doi.org/}%
\providecommand \selectlanguage [0]{\@gobble}%
\providecommand \bibinfo  [0]{\@secondoftwo}%
\providecommand \bibfield  [0]{\@secondoftwo}%
\providecommand \translation [1]{[#1]}%
\providecommand \BibitemOpen [0]{}%
\providecommand \bibitemStop [0]{}%
\providecommand \bibitemNoStop [0]{.\EOS\space}%
\providecommand \EOS [0]{\spacefactor3000\relax}%
\providecommand \BibitemShut  [1]{\csname bibitem#1\endcsname}%
\let\auto@bib@innerbib\@empty
\bibitem [{\citenamefont {Wager}(1911)}]{ref2}%
  \BibitemOpen
  \bibfield  {author} {\bibinfo {author} {\bibfnamefont {H.~W.~T.}\ \bibnamefont {Wager}},\ }\bibfield  {title} {\enquote {\bibinfo {title} {On the effect of gravity upon the movements and aggregation of euglena viridis, ehrb., and other micro-organisms},}\ }\href@noop {} {\bibfield  {journal} {\bibinfo  {journal} {Phil. Trans. R. Soc. Lond. B}\ }\textbf {\bibinfo {volume} {201}},\ \bibinfo {pages} {333--390} (\bibinfo {year} {1911})}\BibitemShut {NoStop}%
\bibitem [{\citenamefont {Platt}(1961)}]{ref3}%
  \BibitemOpen
  \bibfield  {author} {\bibinfo {author} {\bibfnamefont {J.~R.}\ \bibnamefont {Platt}},\ }\bibfield  {title} {\enquote {\bibinfo {title} {``bioconvection patterns" in cultures of free-swimming organisms},}\ }\href@noop {} {\bibfield  {journal} {\bibinfo  {journal} {Science}\ }\textbf {\bibinfo {volume} {133}},\ \bibinfo {pages} {1766--1767} (\bibinfo {year} {1961})}\BibitemShut {NoStop}%
\bibitem [{\citenamefont {Pedley}\ and\ \citenamefont {Kessler}(1992)}]{ref1}%
  \BibitemOpen
  \bibfield  {author} {\bibinfo {author} {\bibfnamefont {T.~J.}\ \bibnamefont {Pedley}}\ and\ \bibinfo {author} {\bibfnamefont {J.~O.}\ \bibnamefont {Kessler}},\ }\bibfield  {title} {\enquote {\bibinfo {title} {Hydrodynamic phenomena in suspensions of swimming micro-organisms},}\ }\href@noop {} {\bibfield  {journal} {\bibinfo  {journal} {Annu. Rev. Fluid Mech.}\ }\textbf {\bibinfo {volume} {24}},\ \bibinfo {pages} {313--358} (\bibinfo {year} {1992})}\BibitemShut {NoStop}%
\bibitem [{\citenamefont {Brinkmann}(1968)}]{ref6}%
  \BibitemOpen
  \bibfield  {author} {\bibinfo {author} {\bibfnamefont {K.}~\bibnamefont {Brinkmann}},\ }\bibfield  {title} {\enquote {\bibinfo {title} {An phasengrenzen induzierte ein und zweidimensionale kristallmuster in kulturen von euglena gracilis},}\ }\href@noop {} {\bibfield  {journal} {\bibinfo  {journal} {Z. Pflanzen Physiol.}\ }\textbf {\bibinfo {volume} {59}},\ \bibinfo {pages} {364--376} (\bibinfo {year} {1968})}\BibitemShut {NoStop}%
\bibitem [{\citenamefont {Kessler}(985b)}]{ref5}%
  \BibitemOpen
  \bibfield  {author} {\bibinfo {author} {\bibfnamefont {J.~O.}\ \bibnamefont {Kessler}},\ }\bibfield  {title} {\enquote {\bibinfo {title} {Co-operative and concentrative phenomena of swimming microorganisms},}\ }\href@noop {} {\bibfield  {journal} {\bibinfo  {journal} {Contemp. Phys.}\ }\textbf {\bibinfo {volume} {26}},\ \bibinfo {pages} {147--166} (\bibinfo {year} {1985b})}\BibitemShut {NoStop}%
\bibitem [{\citenamefont {Nultsch}\ and\ \citenamefont {Hoff}(1993)}]{ref4}%
  \BibitemOpen
  \bibfield  {author} {\bibinfo {author} {\bibfnamefont {W.}~\bibnamefont {Nultsch}}\ and\ \bibinfo {author} {\bibfnamefont {E.}~\bibnamefont {Hoff}},\ }\bibfield  {title} {\enquote {\bibinfo {title} {Investigations on pattern formatin in euglenae},}\ }\href@noop {} {\bibfield  {journal} {\bibinfo  {journal} {Arch. Protistenk}\ }\textbf {\bibinfo {volume} {115}},\ \bibinfo {pages} {336--352} (\bibinfo {year} {1993})}\BibitemShut {NoStop}%
\bibitem [{\citenamefont {Williams}\ and\ \citenamefont {Bees}(2011)}]{ref7}%
  \BibitemOpen
  \bibfield  {author} {\bibinfo {author} {\bibfnamefont {C.~R.}\ \bibnamefont {Williams}}\ and\ \bibinfo {author} {\bibfnamefont {M.~A.}\ \bibnamefont {Bees}},\ }\bibfield  {title} {\enquote {\bibinfo {title} {A tale of three taxes: Photo-gyro-gravitactic bioconvection},}\ }\href@noop {} {\bibfield  {journal} {\bibinfo  {journal} {J. Exp. Biol.}\ }\textbf {\bibinfo {volume} {214}},\ \bibinfo {pages} {2398--2408} (\bibinfo {year} {2011})}\BibitemShut {NoStop}%
\bibitem [{\citenamefont {Häder}(1987)}]{ref8}%
  \BibitemOpen
  \bibfield  {author} {\bibinfo {author} {\bibfnamefont {D.~P.}\ \bibnamefont {Häder}},\ }\bibfield  {title} {\enquote {\bibinfo {title} {Polarotaxis, gravitaxis and vertical phototaxis in the green flagellate, euglena gracilis},}\ }\href@noop {} {\bibfield  {journal} {\bibinfo  {journal} {Arch. Microbiol.}\ }\textbf {\bibinfo {volume} {147}},\ \bibinfo {pages} {179--183} (\bibinfo {year} {1987})}\BibitemShut {NoStop}%
\bibitem [{\citenamefont {Straughan}(1993)}]{ref10}%
  \BibitemOpen
  \bibfield  {author} {\bibinfo {author} {\bibfnamefont {B.}~\bibnamefont {Straughan}},\ }\enquote {\bibinfo {title} {\textit{Mathematical aspects of penetrative convection}},}\ \ (\bibinfo  {publisher} {Longman Scientific},\ \bibinfo {address} {New York},\ \bibinfo {year} {1993})\BibitemShut {NoStop}%
\bibitem [{\citenamefont {Kuznetsov}(2005)}]{kuznetsov2005onset}%
  \BibitemOpen
  \bibfield  {author} {\bibinfo {author} {\bibfnamefont {A.}~\bibnamefont {Kuznetsov}},\ }\bibfield  {title} {\enquote {\bibinfo {title} {The onset of bioconvection in a suspension of gyrotactic microorganisms in a fluid layer of finite depth heated from below},}\ }\href@noop {} {\bibfield  {journal} {\bibinfo  {journal} {International Communications in Heat and Mass Transfer}\ }\textbf {\bibinfo {volume} {32}},\ \bibinfo {pages} {574--582} (\bibinfo {year} {2005})}\BibitemShut {NoStop}%
\bibitem [{\citenamefont {Kuznetsov}(2011)}]{kuznetsov2011non}%
  \BibitemOpen
  \bibfield  {author} {\bibinfo {author} {\bibfnamefont {A.}~\bibnamefont {Kuznetsov}},\ }\bibfield  {title} {\enquote {\bibinfo {title} {Non-oscillatory and oscillatory nanofluid bio-thermal convection in a horizontal layer of finite depth},}\ }\href@noop {} {\bibfield  {journal} {\bibinfo  {journal} {European Journal of Mechanics-B/Fluids}\ }\textbf {\bibinfo {volume} {30}},\ \bibinfo {pages} {156--165} (\bibinfo {year} {2011})}\BibitemShut {NoStop}%
\bibitem [{\citenamefont {Zhao}, \citenamefont {Xiao},\ and\ \citenamefont {Wang}(2018)}]{zhao2018linear}%
  \BibitemOpen
  \bibfield  {author} {\bibinfo {author} {\bibfnamefont {M.}~\bibnamefont {Zhao}}, \bibinfo {author} {\bibfnamefont {Y.}~\bibnamefont {Xiao}}, \ and\ \bibinfo {author} {\bibfnamefont {S.}~\bibnamefont {Wang}},\ }\bibfield  {title} {\enquote {\bibinfo {title} {Linear stability of thermal-bioconvection in a suspension of gyrotactic micro-organisms},}\ }\href@noop {} {\bibfield  {journal} {\bibinfo  {journal} {International Journal of Heat and Mass Transfer}\ }\textbf {\bibinfo {volume} {126}},\ \bibinfo {pages} {95--102} (\bibinfo {year} {2018})}\BibitemShut {NoStop}%
\bibitem [{\citenamefont {Zhao}\ \emph {et~al.}(2019)\citenamefont {Zhao}, \citenamefont {Wang}, \citenamefont {Wang},\ and\ \citenamefont {Mahabaleshwar}}]{zhao2019darcy}%
  \BibitemOpen
  \bibfield  {author} {\bibinfo {author} {\bibfnamefont {M.}~\bibnamefont {Zhao}}, \bibinfo {author} {\bibfnamefont {S.}~\bibnamefont {Wang}}, \bibinfo {author} {\bibfnamefont {H.}~\bibnamefont {Wang}}, \ and\ \bibinfo {author} {\bibfnamefont {U.}~\bibnamefont {Mahabaleshwar}},\ }\bibfield  {title} {\enquote {\bibinfo {title} {Darcy--brinkman bio-thermal convection in a suspension of gyrotactic microorganisms in a porous medium},}\ }\href@noop {} {\bibfield  {journal} {\bibinfo  {journal} {Neural Computing and Applications}\ }\textbf {\bibinfo {volume} {31}},\ \bibinfo {pages} {1061--1067} (\bibinfo {year} {2019})}\BibitemShut {NoStop}%
\bibitem [{\citenamefont {Balla}\ \emph {et~al.}(2020)\citenamefont {Balla}, \citenamefont {Ramesh}, \citenamefont {Kishan}, \citenamefont {Rashad},\ and\ \citenamefont {Abdelrahman}}]{balla2020bioconvection}%
  \BibitemOpen
  \bibfield  {author} {\bibinfo {author} {\bibfnamefont {C.~S.}\ \bibnamefont {Balla}}, \bibinfo {author} {\bibfnamefont {A.}~\bibnamefont {Ramesh}}, \bibinfo {author} {\bibfnamefont {N.}~\bibnamefont {Kishan}}, \bibinfo {author} {\bibfnamefont {A.}~\bibnamefont {Rashad}}, \ and\ \bibinfo {author} {\bibfnamefont {Z.}~\bibnamefont {Abdelrahman}},\ }\bibfield  {title} {\enquote {\bibinfo {title} {Bioconvection in oxytactic microorganism-saturated porous square enclosure with thermal radiation impact},}\ }\href@noop {} {\bibfield  {journal} {\bibinfo  {journal} {Journal of Thermal Analysis and Calorimetry}\ }\textbf {\bibinfo {volume} {140}},\ \bibinfo {pages} {2387--2395} (\bibinfo {year} {2020})}\BibitemShut {NoStop}%
\bibitem [{\citenamefont {Hussain}, \citenamefont {Raizah},\ and\ \citenamefont {Aly}(2022)}]{hussain2022thermal}%
  \BibitemOpen
  \bibfield  {author} {\bibinfo {author} {\bibfnamefont {S.}~\bibnamefont {Hussain}}, \bibinfo {author} {\bibfnamefont {Z.}~\bibnamefont {Raizah}}, \ and\ \bibinfo {author} {\bibfnamefont {A.~M.}\ \bibnamefont {Aly}},\ }\bibfield  {title} {\enquote {\bibinfo {title} {Thermal radiation impact on bioconvection flow of nano-enhanced phase change materials and oxytactic microorganisms inside a vertical wavy porous cavity},}\ }\href@noop {} {\bibfield  {journal} {\bibinfo  {journal} {International Communications in Heat and Mass Transfer}\ }\textbf {\bibinfo {volume} {139}},\ \bibinfo {pages} {106454} (\bibinfo {year} {2022})}\BibitemShut {NoStop}%
\bibitem [{\citenamefont {Vincent}\ and\ \citenamefont {Hill}(1996)}]{ref9}%
  \BibitemOpen
  \bibfield  {author} {\bibinfo {author} {\bibfnamefont {R.~V.}\ \bibnamefont {Vincent}}\ and\ \bibinfo {author} {\bibfnamefont {N.~A.}\ \bibnamefont {Hill}},\ }\bibfield  {title} {\enquote {\bibinfo {title} {Bioconvection in a suspension of phototactic algae},}\ }\href@noop {} {\bibfield  {journal} {\bibinfo  {journal} {J. Fluid Mech.}\ }\textbf {\bibinfo {volume} {327}},\ \bibinfo {pages} {343--371} (\bibinfo {year} {1996})}\BibitemShut {NoStop}%
\bibitem [{\citenamefont {Ghorai}, \citenamefont {Panda},\ and\ \citenamefont {Hill}(2010)}]{ref14}%
  \BibitemOpen
  \bibfield  {author} {\bibinfo {author} {\bibfnamefont {S.}~\bibnamefont {Ghorai}}, \bibinfo {author} {\bibfnamefont {M.~K.}\ \bibnamefont {Panda}}, \ and\ \bibinfo {author} {\bibfnamefont {N.~A.}\ \bibnamefont {Hill}},\ }\bibfield  {title} {\enquote {\bibinfo {title} {Bioconvection in a suspension of isotropically scattering phototactic algae},}\ }\href@noop {} {\bibfield  {journal} {\bibinfo  {journal} {Phys. Fluids}\ }\textbf {\bibinfo {volume} {22}},\ \bibinfo {pages} {071901} (\bibinfo {year} {2010})}\BibitemShut {NoStop}%
\bibitem [{\citenamefont {Kumar}(2023{\natexlab{a}})}]{kumar2023}%
  \BibitemOpen
  \bibfield  {author} {\bibinfo {author} {\bibfnamefont {S.}~\bibnamefont {Kumar}},\ }\bibfield  {title} {\enquote {\bibinfo {title} {Isotropic scattering with a rigid upper surface at the onset of phototactic bioconvection},}\ }\href@noop {} {\bibfield  {journal} {\bibinfo  {journal} {Physics of Fluids}\ }\textbf {\bibinfo {volume} {35}},\ \bibinfo {pages} {024106} (\bibinfo {year} {2023}{\natexlab{a}})}\BibitemShut {NoStop}%
\bibitem [{\citenamefont {Ghorai}\ and\ \citenamefont {Panda}(2013)}]{ref15}%
  \BibitemOpen
  \bibfield  {author} {\bibinfo {author} {\bibfnamefont {S.}~\bibnamefont {Ghorai}}\ and\ \bibinfo {author} {\bibfnamefont {M.~K.}\ \bibnamefont {Panda}},\ }\bibfield  {title} {\enquote {\bibinfo {title} {Bioconvection in an anisotropic scattering suspension of phototactic algae},}\ }\href@noop {} {\bibfield  {journal} {\bibinfo  {journal} {Eur. J. Mech.-B/Fluids}\ }\textbf {\bibinfo {volume} {41}},\ \bibinfo {pages} {81--93} (\bibinfo {year} {2013})}\BibitemShut {NoStop}%
\bibitem [{\citenamefont {Panda}\ \emph {et~al.}(2016)\citenamefont {Panda}, \citenamefont {Singh}, \citenamefont {Mishra},\ and\ \citenamefont {Mohanty}}]{ref16}%
  \BibitemOpen
  \bibfield  {author} {\bibinfo {author} {\bibfnamefont {M.~K.}\ \bibnamefont {Panda}}, \bibinfo {author} {\bibfnamefont {R.}~\bibnamefont {Singh}}, \bibinfo {author} {\bibfnamefont {A.~C.}\ \bibnamefont {Mishra}}, \ and\ \bibinfo {author} {\bibfnamefont {S.~K.}\ \bibnamefont {Mohanty}},\ }\bibfield  {title} {\enquote {\bibinfo {title} {Effects of both diffuse and collimated incident radiation on phototactic bioconvection},}\ }\href@noop {} {\bibfield  {journal} {\bibinfo  {journal} {Phys. Fluids}\ }\textbf {\bibinfo {volume} {28}},\ \bibinfo {pages} {124104} (\bibinfo {year} {2016})}\BibitemShut {NoStop}%
\bibitem [{\citenamefont {Panda}(2020)}]{panda2020effects}%
  \BibitemOpen
  \bibfield  {author} {\bibinfo {author} {\bibfnamefont {M.}~\bibnamefont {Panda}},\ }\bibfield  {title} {\enquote {\bibinfo {title} {Effects of anisotropic scattering on the onset of phototactic bioconvection with diffuse and collimated irradiation},}\ }\href@noop {} {\bibfield  {journal} {\bibinfo  {journal} {Physics of Fluids}\ }\textbf {\bibinfo {volume} {32}} (\bibinfo {year} {2020})}\BibitemShut {NoStop}%
\bibitem [{\citenamefont {Panda}, \citenamefont {Sharma},\ and\ \citenamefont {Kumar}(2022)}]{ref17}%
  \BibitemOpen
  \bibfield  {author} {\bibinfo {author} {\bibfnamefont {M.~K.}\ \bibnamefont {Panda}}, \bibinfo {author} {\bibfnamefont {P.}~\bibnamefont {Sharma}}, \ and\ \bibinfo {author} {\bibfnamefont {S.}~\bibnamefont {Kumar}},\ }\bibfield  {title} {\enquote {\bibinfo {title} {Effect of oblique irradiation on the onset of phototactic bioconvection},}\ }\href@noop {} {\bibfield  {journal} {\bibinfo  {journal} {Phys. Fluids}\ }\textbf {\bibinfo {volume} {34}},\ \bibinfo {pages} {024108} (\bibinfo {year} {2022})}\BibitemShut {NoStop}%
\bibitem [{\citenamefont {Kumar}(2022)}]{ref18}%
  \BibitemOpen
  \bibfield  {author} {\bibinfo {author} {\bibfnamefont {S.}~\bibnamefont {Kumar}},\ }\bibfield  {title} {\enquote {\bibinfo {title} {Phototactic isotropic scattering bioconvection with oblique irradiation},}\ }\href@noop {} {\bibfield  {journal} {\bibinfo  {journal} {Phys. Fluids}\ }\textbf {\bibinfo {volume} {34}},\ \bibinfo {pages} {114125} (\bibinfo {year} {2022})}\BibitemShut {NoStop}%
\bibitem [{\citenamefont {Kumar}(2023{\natexlab{b}})}]{kumar2023effect}%
  \BibitemOpen
  \bibfield  {author} {\bibinfo {author} {\bibfnamefont {S.}~\bibnamefont {Kumar}},\ }\bibfield  {title} {\enquote {\bibinfo {title} {Effect of rotation on the suspension of phototactic bioconvection},}\ }\href@noop {} {\bibfield  {journal} {\bibinfo  {journal} {Physics of Fluids}\ }\textbf {\bibinfo {volume} {35}} (\bibinfo {year} {2023}{\natexlab{b}})}\BibitemShut {NoStop}%
\bibitem [{\citenamefont {Ghorai}\ and\ \citenamefont {Hill}(2005)}]{ref13}%
  \BibitemOpen
  \bibfield  {author} {\bibinfo {author} {\bibfnamefont {S.}~\bibnamefont {Ghorai}}\ and\ \bibinfo {author} {\bibfnamefont {N.~A.}\ \bibnamefont {Hill}},\ }\bibfield  {title} {\enquote {\bibinfo {title} {Penetrative phototactic bioconvection},}\ }\href@noop {} {\bibfield  {journal} {\bibinfo  {journal} {Phys. Fluids}\ }\textbf {\bibinfo {volume} {17}},\ \bibinfo {pages} {074101} (\bibinfo {year} {2005})}\BibitemShut {NoStop}%
\bibitem [{\citenamefont {Modest}(2003)}]{ref-modest}%
  \BibitemOpen
  \bibfield  {author} {\bibinfo {author} {\bibfnamefont {M.~F.}\ \bibnamefont {Modest}},\ }\enquote {\bibinfo {title} {\textit{Radiative Heat Transfer}},}\ \ (\bibinfo  {publisher} {Academic Press},\ \bibinfo {address} {New York},\ \bibinfo {year} {2003})\ \bibinfo {edition} {2nd}\ ed.\BibitemShut {Stop}%
\bibitem [{\citenamefont {Chandrasekhar}(1960)}]{ref-chand}%
  \BibitemOpen
  \bibfield  {author} {\bibinfo {author} {\bibfnamefont {S.}~\bibnamefont {Chandrasekhar}},\ }\enquote {\bibinfo {title} {\textit{Radiative Transfer}},}\ \ (\bibinfo  {publisher} {Dover},\ \bibinfo {address} {New York},\ \bibinfo {year} {1960})\BibitemShut {NoStop}%
\bibitem [{\citenamefont {Shampine}, \citenamefont {Gladwell},\ and\ \citenamefont {Thompson}(2003)}]{shampine2003solving}%
  \BibitemOpen
  \bibfield  {author} {\bibinfo {author} {\bibfnamefont {L.~F.}\ \bibnamefont {Shampine}}, \bibinfo {author} {\bibfnamefont {I.}~\bibnamefont {Gladwell}}, \ and\ \bibinfo {author} {\bibfnamefont {S.}~\bibnamefont {Thompson}},\ }\href@noop {} {\emph {\bibinfo {title} {Solving ODEs with matlab}}}\ (\bibinfo  {publisher} {Cambridge university press},\ \bibinfo {year} {2003})\BibitemShut {NoStop}%
\bibitem [{\citenamefont {Chandrasekhar}(1961)}]{ref12}%
  \BibitemOpen
  \bibfield  {author} {\bibinfo {author} {\bibfnamefont {S.}~\bibnamefont {Chandrasekhar}},\ }\enquote {\bibinfo {title} {\textit{Hydrodynamic and hydromagnetic stability}},}\ \ (\bibinfo  {publisher} {Oxford University Press},\ \bibinfo {year} {1961})\BibitemShut {NoStop}%
\end{thebibliography}%


\end{document}